\documentclass[12pt, twoside, leqno]{article}

\usepackage{amsmath,amsthm}
\usepackage{amssymb}

\usepackage{enumerate}

\newtheorem{thm}{Theorem}[section]
\newtheorem{cor}[thm]{Corollary}
\newtheorem{lem}[thm]{Lemma}
\newtheorem{alg}[thm]{Algorithm}
\newtheorem{prob}[thm]{Problem}

\theoremstyle{definition}
\newtheorem{defin}[thm]{Definition}
\newtheorem{rem}[thm]{Remark}
\newtheorem{exa}[thm]{Example}

\newcommand{\ZZ}{\mathbb{Z}}
\newcommand{\NN}{\mathbb{N}}

\newcommand{\PP}{\mathbb{P}}
\DeclareMathOperator{\lc}{lc}

\begin{document}


\baselineskip=17pt


\title{Digit Polynomials and their application to integer factorization}

\author{Markus Hittmeir\\
University of Salzburg\\
Mathematics Department}
\date{}

\maketitle


\renewcommand{\thefootnote}{}

\footnote{The author is supported by the Austrian Science Fund (FWF): Project F5504-N26.}

\footnote{Address: Hellbrunnerstra{\ss}e 34, A-5020 Salzburg. E-Mail: markus.hittmeir@sbg.ac.at}

\footnote{2010 \emph{Mathematics Subject Classification}: Primary 11A51; Secondary 11A41.}

\footnote{\emph{Key words and phrases}: Factorization, Primality, Primes.}

\renewcommand{\thefootnote}{\arabic{footnote}}
\setcounter{footnote}{0}


\vspace{-2.5cm}

\begin{abstract}
This paper presents the concept of \emph{digit polynomials}, which leads to a deterministic and unconditional integer factorization algorithm with the runtime complexity $\mathcal{O}(N^{1/4+\epsilon})$. Strassen's well known factoring approach is a special case of our method. We will also consider a possibility to improve upon the complexity bound. 
\end{abstract}

\section{Introduction}
We consider the problem of computing the prime factorization of a given natural number $N$. Currently, the best publicly known deterministic and unconditional factorization algorithms all have a runtime complexity of the form $\mathcal{O}(N^{1/4+\epsilon})$ \cite[p.240]{Wag}. A method which achieves this complexity is the approach of Strassen \cite{Str}, based on the idea to compute parts of $\lfloor N^{1/2}\rfloor!$ to find a nontrivial factor of $N$. A recent improvement of the logarithmic factor in the complexity bound can be found in \cite{CosHar}. For a general overview, the reader may consult \cite{Pom}.

In this paper we present a method based on products of certain polynomials. The main idea is to construct polynomials $g\in\ZZ[X]$ such that as many integers $x$, $0\leq x\leq N-1$, as possible satisfy 
\[1<\gcd(g(x),N)<N.\]
Several $b$-adic representations of $N$ are used in Theorem \ref{d}, which yields a method to construct such a polynomial of degree $d$ with complexity $\mathcal{O}(d^{1+\epsilon})$. In the factorization algorithm we will not only make use of the \emph{cardinality}, but also of the \emph{position} of those $x$ with the property above.

Our deterministic method is not appropriate for factorizing large numbers. In practice, probabilistic algorithms with much lower complexity are used for this task (See \cite{Rie} and \cite{CraPom}).

\section{Basic Ideas}

Throughout this paper, $\PP$ denotes the set of primes. We call a natural number semiprime if and only if it is the product of two distinct primes. Let $n\in\NN$. We denote the complete residue system $\{0,...,n-1\}$ modulo $n$ by $Z_n$ and the residue class ring $\ZZ/n\ZZ$ by $\ZZ_n$. For $f\in\ZZ[X]$, we write the leading coefficient of $f$ as $\lc f$. Until further notice, let $N\in\NN$ be fixed. 

\defin{Let $b\in\ZZ$. We denote the set of polynomials $f\in\ZZ[X]$ with the property $f(b)=N$ by $\mathcal{D}_{N,b}$. The elements of $\mathcal{D}_{N,b}$ are called \emph{digit polynomials of $N$ to base $b$}.}

\defin{Let $b\in\NN$, $b\geq 2$. Let $N=\sum_{i\geq 0} n_ib^{i}$ be the unique $b$-adic representation of $N$ with digits $n_i\in\{0,...,b-1\}$. Define
\[
  P_b:=\sum_{i\geq 0}n_iX^{i}\in\ZZ[X].
\]
We call $P_b$ the \emph{b-adic digit polynomial of $N$}. Clearly, we have $P_b\in\mathcal{D}_{N,b}$.}

{\lem{Let $b\in\ZZ$ and $f\in\mathcal{D}_{N,b}$. Then, for every $x\in\ZZ$, we have $N\equiv f(x)\mod x-b$.\label{a}}}

\begin{proof}
  We know that $b$ is a zero of the polynomial $f-N$, hence $X-b$ divides $f-N$ in $\ZZ[X]$ and the congruence holds for every evaluation.
\end{proof}

{\cor{Let  $b\in\ZZ$ and $f\in\mathcal{D}_{N,b}$. We conclude for every $x\in\ZZ$ that $\gcd(N,x-b)=\gcd(f(x),x-b)$, and that $x-b\mid N$ iff $x-b\mid f(x)$.}}

{\lem{Let $u$ and $v$ be nontrivial and coprime divisors of $N$. Let $b\in\ZZ$ and $f\in\mathcal{D}_{N,b}$ such that
\begin{enumerate}
\item{$\gcd(\lc f,N)=1$ and}
\item{$d:=\deg f$ is smaller than the largest prime factor of $v$.}
\end{enumerate} 
Then there exists $x\in\ZZ$ with $u\mid f(x)$ and $v\nmid f(x)$.\label{b}}}

\begin{proof}
 Let $y\in\ZZ$ be arbitrary. Let $x\in\ZZ$ with $uy=x-b$. From Lemma \ref{a} we derive $u\mid f(x)$, hence $u\mid f(uy+b)$ for any $y\in\ZZ$. We have to show that there exists $y\in\ZZ$ with $v\nmid f(uy+b)$.
 
 Assume to the contrary that $f(uy+b)\equiv 0 \mod v$ for all $y\in\ZZ$. Write $f(uy+b)$ as $f(b)+u\cdot g(y)$ for $g\in\ZZ[X]$. It is easy to verify that $\deg g=d$ and $\lc g=u^{d-1}\lc f$. Let $p$ be the largest prime factor of $v$. Then, for every $y\in\ZZ$, it follows that
\[
   f(uy+b)=u\cdot g(y)+f(b)\equiv u\cdot g(y)\equiv 0 \mod p.
\]
The fact $p\nmid u$ implies $g(y)\equiv 0 \mod p$ for every $y\in\ZZ$. But, since $\gcd(\lc f,N)=1$, we get $p\nmid lc(g)$. Therefore, $g$ is of degree $d$ in $\ZZ_p[X]$ and, for this reason, has at most $d$ zeros in $\ZZ_p[X]$. From $d<p$ the contradiction follows.
\end{proof}

In the proof of the preceding lemma we have seen that, if $N$ is a composite number and if $f\in\mathcal{D}_{N,b}$ is chosen with appropriate degree, we get various integers $x\in Z_N$ such that $1<\gcd(f(x),N)<N$.

\defin{Let $g\in\ZZ[X]$. An element $x\in Z_N$ is called \emph{suitable for $g$}, if and only if $1<\gcd(g(x),N)<N$. We also define
\[
  \nu(g):=\#\{x\in Z_N:\text{$x$ is suitable for $g$.}\}.
\]
}
If we multiply two polynomials $f,g\in\ZZ[X]$, it may happen that $x\in Z_N$ is suitable for $f$ and for $g$, but not for $f\cdot g$.

\defin{Let $d\in\NN$ and $f_i\in\ZZ[X]$, $1\leq i \leq d$.  An element $x\in Z_N$ \emph{vanishes in $g:=\prod_{i=1}^d f_i$}, if and only if $\gcd(g(x),N)=N$ and there is at least one $i$ such that $x$ is suitable for $f_i$.}

{\thm{Let $N\in\NN$ be a semiprime number with the prime factors $p$ and $q$ and assume $p<q$. Let $f\in\ZZ[X]$ and $d:=\deg f$. Let $n$ be the number of distinct zeros of $f$ modulo $p$ and $m$ be the number of distinct zeros of $f$ modulo $q$. Then:

\begin{enumerate}
  \item{$\nu(f)=mp+nq-2nm$.\label{thm1:1}}
  \item{Let $f\neq 0$ in $\ZZ_p[X]$ and in $\ZZ_q[X]$. If $d<p/2$, then $\nu (f)\leq dp+dq-2d^2$.\label{thm1:2}}
\end{enumerate} \label{c}}}

\begin{proof}
 For \ref{thm1:1}: Let $x\in Z_N$ be suitable for $f$. Then $x$ is a zero of $f$ either modulo $p$ or modulo $q$. Let $\alpha_1,...,\alpha_n$ be the distinct zeros of $f$ modulo $p$ and $\beta_1,...,\beta_m$ be the distinct zeros of $f$ modulo $q$. For $i=1,...,n$ and $j=1,...,m$ we consider
  \begin{align*}
    &py+\alpha_i, \text{ for } y=0,...,q-1,\\
    &qy+\beta_j, \text{ for } y=0,...,p-1.
  \end{align*}
Every $x$ which is suitable for $f$ is of that form, and these are a priori $mp+nq$ values in $Z_N$. But some of them might be equal. First, we show that the values of the form $py+\alpha_i$ are distinct modulo $N$. We assume that there are $y_1,y_2\in Z_q$ with $py_1+\alpha_i \equiv py_2+\alpha_k \mod N$ for some $i,k\in \{1,...,n\}$. For $i\neq k$ this is not possible, because we get $\alpha_i\equiv \alpha_j \mod p$, which contradicts the assumption that the zeros are distinct modulo $p$. For $i=k$, it follows that $y_1\equiv y_2 \mod q$. Hence, the congruence only holds if we compare the value $py_1+\alpha_i$ with itself. For this reason, all these values are distinct. By similar arguments, one can show that this also holds for the values of the form $qy+\beta_j$.

Next, we consider the case that some value of the form $py+\alpha_i$ is congruent to some value of the form $qy+\beta_j$. Then this value is a zero of $f$ modulo $N$. By the Chinese Remainder Theorem, one can easily verify that $f$ must have exactly $nm$ distinct zeros modulo $N$. Since any zero $z$ of $f$ modulo $N$ is also a zero of $f$ modulo $p$ and modulo $q$, we can write $z=py_1+\alpha_i=qy_2+\beta_j$ for some $y_1,y_2$ and $i,j$. Hence, at every zero of $f$ modulo $N$ exactly two equal values of our list above coincide. The other values all satisfy $1<\gcd(f(x),N)<N$. Therefore, we get $\nu(f)=mp+nq-2nm$.\\

For \ref{thm1:2}: Consider $h=-2XY+Xq+Yp\in\ZZ[X,Y]$. Since $f$ has at most $d$ distinct zeros modulo $p$ and modulo $q$, we want to maximize this function for $(x,y)\in [0,d]^2$. We get
\[
  h_X(x,y)=-2y+q \text{ and } h_Y(x,y)=-2x+p
\]
as partial derivatives. Hence, the only critical point is $(x,y)=(p/2,q/2)$. But this point is not in $[0,d]^2$, so we consider $h$ on the boundary and get $g_1(x)=xq$, $g_2(x)=xp$, $g_3(x)=x(q-2d)+dp$ and $g_4(x)=x(p-2d)+dq$, for $x\in [0,c]$. Since $q>2d$ and $p>2d$, the maximum is 
\[
g_3(d)=g_4(d)=dp+dq-2d^2.
\].
\end{proof}

For any polynomial $f\in\ZZ[X]$ with appropriate degree $d$, there are at most $dp+dq-2d^2$ integers which are suitable for $f$. We are interested in \emph{efficient methods} to construct polynomials which are best possible in this sense. The following theorem yields a method with runtime complexity of the form $\mathcal{O}(d^{1+\epsilon})$. We will use this idea in the factorization algorithm in Section 3. Therefore, details will be explained in the proof of Theorem \ref{n}.

{\thm{Let $N\in\NN$ be semiprime with the prime factors $p$ and $q$. Let $d\in\NN$ and $b_i\in\ZZ$, $1\leq i\leq d$. Let $f_i\in\mathcal{D}_{N,b_i}$ such that $\deg f_i=1$ and write $f_i=l_iX+c_i$ for every $i$. If $\gcd(c_i,N)=1$ for every $i$ and also $\gcd(b_j-b_k,N)=1$ for every choice of $j,k\in\{1,...,d\}$, $j\neq k$, then
\[
  \nu\Big{(}\prod_{i=1}^d f_i\Big{)}=dp+dq-2d^2.
\]
}\label{d}}

\begin{proof}
 For $1\leq i\leq d$, consider $f_i$. Since $\gcd(c_i,N)=1$, $f_i\neq 0$ as polynomial in $\ZZ_p[X]$ and in $\ZZ_q[X]$. Therefore, $b_i$ is the only zero of $f_i$ modulo $p$ and modulo $q$.
  
  Now consider $g:=\prod_{i=1}^d f_i$. Obviously, every $b_i$, $1\leq i\leq d$, is a zero of $g$ modulo $p$ as well as modulo $q$. Since $\gcd(b_j-b_k,N)=1$ for every choice of $j,k\in\{1,...,d\}$, $j\neq k$, these zeros are distinct. For this reason, $g$ has $d$ distinct zeros modulo $p$ and $d$ distinct zeros modulo $q$. Now we apply Theorem \ref{c}.
\end{proof}

\rem{For every polynomial $f_i$ in the theorem above, there are $p+q-2$ integers which are suitable for $f_i$. But, if we multiply all these polynomials, we do not get $d(p+q-2)$ suitable integers for the product $g$. It is easy to see that there are $4\cdot {d \choose 2}$ integers vanishing in $g$. We get
\[
  \nu\Big{(}\prod_{i=1}^df_i\Big{)}=dp+dq-2d^2=d(p+q-2)-4\cdot {d \choose 2}.
\]}

We will now prove a result to ensure the maximum possible number of suitable integers for the product of digit polynomials of degree $2$, which may be compared to the result in Theorem \ref{d}. We will see that the $b$-adic digit polynomials are especially useful in this case, not only because it is easy to compute them, but also because of their uniqueness and the special way they are constructed.

{\thm{Let $N\in\NN$ be semiprime with prime factors $p$ and $q$. Let $d\in\NN$ and $b_i\in\ZZ$, $1\leq i\leq d$. Let $f_i\in\mathcal{D}_{N,b_i}$ such that $\deg f_i=2$ and write $f_i=n_{2,i}X^2+n_{1,i}X+n_{0,i}$ for every $i$.

If $\gcd(n_{2,i}\cdot b_i,N)=1$ for every $i$ and if for $b_{d+i}:=n_{0,i}\cdot n_{2,i}^{-1}\cdot b_i^{-1} \mod N$, $1\leq i\leq d$, we have $\gcd(b_j-b_k,N)=1$ for every choice of $j,k\in\{1,...,2d\}$, $j\neq k$, then
\[
\nu\Big{(}\prod_{i=1}^d f_i\Big{)}=2dp+2dq-8d^2.
\]
}\label{g}}

\begin{proof}
For $1\leq i\leq d$, consider $f_i$. Since $\gcd(n_{2,i},N)=1$, $f_i$ is a polynomial of degree $2$ modulo $p$. Therefore $f_i$ has at most two zeros modulo $p$. One of them is $b_i$. But since $\ZZ_p$ is a field, there has to be another zero modulo $p$. We know from Vieta's Theorem that this zero has to be the solution of 
\[
n_{2,i}b_i\cdot x\equiv n_{0,i}\mod p.
\]
Since $b_{d+i}\equiv n_{0,i}\cdot n_{2,i}^{-1}\cdot b_i^{-1} \mod p$, $b_{d+i}$ is this zero of $f_i$ modulo $p$. With similar arguments, one can also show that $b_i$ and $b_{d+i}$ are the zeros of $f_i$ modulo $q$.

Now consider $g:=\prod_{i=1}^d f_i$. Obviously, every $b_i$, $1\leq i \leq 2d$ is a zero of $g$ modulo $p$ as well as modulo $q$. Since $\gcd(b_j-b_k,N)=1$ for every choice of $j,k\in\{1,...,2d\}$, $j\neq k$, these zeros are distinct. For this reason, $g$ has $2d$ distinct zeros modulo $p$ and $2d$ distinct zeros modulo $q$. Now we apply Theorem \ref{c}.
\end{proof}

If we set $d=1$ in the theorem above, the following statement is an immediate consequence.

{\cor{Let $N\in\NN$ be semiprime with prime factors $p$ and $q$. Let $b\in\ZZ$ and $f\in\mathcal{D}_{N,b}$ with $\deg f=2$ and $f=n_2X^2+n_1X+n_0$. If $\gcd(n_2\cdot b,N)=1$ and $\gcd(N,n_2b^2-n_0)=1$, then $\nu(f)=2p+2q-8$.}}\\

We want to make Theorem \ref{g} applicable. Hence, we have to find digit polynomials for which the condition of distinct zeros modulo the factors of $N$ can be verified in $\mathcal{O}(d)$ steps. For the linear polynomials in Theorem \ref{d} this is feasible, since we are able to choose appropriate bases, for example consecutive integers. Here, every base $b_i$ we choose comes with a second integer $b_{d+i}$, which we have to control. The subsequent lemma allows to work with digit polynomials of degree $2$ in practice.

{\lem{Let $N\in\NN$, $d\in\NN$ and let $b_i\in\{\lceil N^{1/2}/\sqrt{2}\rceil,...,\lfloor N^{1/2}\rfloor\}$, $1\leq i \leq d$, be coprime to $N$ such that $b_{i+1}=b_i+1$. Set $D:=b_1+\lfloor N/b_d\rfloor.$

If $\gcd(D+z,N)=1$ for every $z\in\{0,...,2d-2\}$ and if the $b$-adic digit polynomials $P_{b_i}$ satisfy $n_{1,i}\leq n_{0,i}+1$ for every $i$, then they also satisfy the conditions in Theorem \ref{g}.
}}

\begin{proof}
Let $i\in\{1,...,d\}$ be arbitrary. It is easy to see that $n_{2,i}=1$ for this choice of bases. Since $\gcd(b_i,N)=1$, the first condition of Theorem \ref{g} is satisfied. Now set $b_{d+i}:=n_{0,i}b_i^{-1}\mod N$. Consider the division with remainder of $N$ with respect to $b_i$ and write $m_ib_i+n_{0,i}=N$. We get $-m_{i}b_i\equiv n_{0,i}\equiv b_{d+i}b_i\mod N$, hence $-m_{i}\equiv b_{d+i}\mod N$. Next, we consider $N=(m_i-1)(b_i+1)+r=m_ib_i+m_i-b_i-1+r$ for some $r\in\ZZ$. Assume that $r\geq b_i+1$. Then it follows that $N\geq m_ib_i+m_i$. But since $b_i\leq \lfloor N^{1/2}\rfloor$, it is easy to see that $m_i\geq b_i$. By $n_{0,i}<b_i$ we conclude
\[
  N\geq m_ib_i+m_i\geq m_ib_i+b_i>m_ib_i+n_{0,i}=N,
\]
hence a contradiction. Now we assume that $r<0$. Then it follows that $N<m_ib_i+m_i-b_i-1$. But this yields that $n_{0,i}+b_i+1<m_i$, and by $N>b_i(n_{0,i}+b_i+1)+n_{0,i}=b_i^2+(n_{0,i}+1)b_i+n_{0,i}$ we conclude $n_{1,i}>n_{0,i}+1$, which contradicts our assumption. As a consequence, we get $0\leq r<b_i+1$. Because of the uniqueness of the division with remainder, there has to be $r=n_{0,i+1}$ and $m_{i+1}=m_i-1$. Altogether we derive
\[
  b_{d+i+1}\equiv -m_{i+1}\equiv -m_i+1\equiv b_{d+i}+1\mod N.
\]
Now assume that there exist $j,k \in\{1,...,d\}$ such that $b_{d+k}\equiv b_j\mod p$. We write $b_j=b_1+m$ for some $m\in\{0,...,d-1\}$ and, as we just have shown, we can write
\[
  b_{d+k}\equiv b_{2d}-l\equiv -m_d-l\equiv -\lfloor N/b_d\rfloor-l\mod p,
\]
for some $l\in\{0,...,d-1\}$. It follows that  $-\lfloor N/b_d\rfloor -l\equiv b_1+m\mod p$. Therefore, we get
\[
   0\equiv b_1+\lfloor N/b_d\rfloor+m+l\equiv D+z\mod p,
\]
for some $z\in\{0,...,2d-2\}$. But this contradicts our assumption. Hence, for every choice of $j,k\in\{1,...,d\}$, the integers $b_{d+j}$ are different from the integers $b_k$ modulo $p$. It is also impossible that there exist $j,k\in\{1,...,d\}$, $j\neq k$ with $b_{d+k}\equiv b_{d+j}\mod p$ or with $b_k\equiv b_j\mod p$, because this would imply $p\leq d$, which as well contradicts the assumption $\gcd(D+z,N)=1$ for $z\in\{0,...,2d-2\}$. By similar arguments, one can show that the zeros are all distinct modulo $q$.
\end{proof}

\section{The Algorithm and its Parameters}
Let $N\in\NN$ be a composite number. Without knowledge of the factorization of $N$, we are able to construct a polynomial $g\in\ZZ[X]$ such that
\[
1<\gcd(g(x),N)<N,
\]
for as many $x\in Z_N$ as possible. The main idea for the algorithm is to find a subset of $Z_N$ containing at least one element which is either suitable for or vanishing in $g$. Let $d\in\NN$. We work with the following parameters. 

\begin{enumerate}
\item{A set $\mathcal{B}:=\{b_n\in Z_N:1\leq n\leq d\}$ of bases for the digit polynomials.}
\item{For every $b\in\mathcal{B}$, we choose exactly one $f_b\in\mathcal{D}_{N,b}$. We denote the set of all these polynomials by $\mathcal{D}(\mathcal{B})$.}
\item{A set $\mathcal{S}:=\{s_n\in Z_N:1\leq n\leq d\}$, containing at least one element suitable for or vanishing in $g:=\prod_{b\in\mathcal{B}} f_b$.}
\end{enumerate}

These three sets determine the following algorithm, and its correctness and runtime depends on finding a good choice for them.

{\alg{Let $N\in\NN$ and the sets $\mathcal{B}=\{b_n\in Z_N:1\leq n\leq d\}$, $\mathcal{D}(\mathcal{B})=\{f_b\in\mathcal{D}_{N,b}:b\in\mathcal{B}\}$ and $\mathcal{S}=\{s_n\in Z_N:1\leq n\leq d\}$ be given, where $d\in\NN$. Set $a_1=1$, $a_2=1$ and take the following steps to factor $N$:

\begin{enumerate}
	\item{For every $b\in\mathcal{B}$, compute $f_b\in\mathcal{D}(\mathcal{B})$. Next, compute the polynomial $g:=\prod_{b\in\mathcal{B}} f_b \mod N$.}
	\item{For every $n\in\{1,...,d\}$, compute $y_n:=g(s_n)\mod N$.}
	\item{Set $j:=a_1$.}
	\item{If $j>d$, print 'Error A'. Otherwise compute $G_j:=\gcd(y_j,N)$. If $G_j=1$, set $a_1=j+1$ and go to Step $3$. If $1<G_j<N$, print $G_j$. We have found a nontrivial factor of $N$ and the algorithm terminates. If $G_{j}=N$, go to Step $5$.}
	\item{Set $i:=a_2$.}
	\item{If $i>d$, print 'Error B'. Otherwise compute $H_i:=\gcd(f_{b_i}(s_j),N)$. If $H_i=1$ or $H_i=N$, set $a_2=i+1$ and go to Step $5$. If $1<H_i<N$, print $H_i$. We have found a nontrivial factor of $N$ and the algorithm terminates.}
	\end{enumerate}
	\label{l}}}
	
We now clarify which conditions are necessary to make the algorithm work. Finding a solution to the following problem is crucial. 
	
{\prob{Let $N\in\NN$ be of unknown factorization. For $d\in\NN$, construct two disjoint sets $\{b_n:1\leq n\leq d\}$ and $\{s_n:1\leq n\leq d\}$ in $Z_N$ with the property that, if $N$ is composite, there must exist $i,j\in\{1,...,d\}$ and a prime factor $p$ of $N$ such that $b_i\equiv s_j\mod p$.
\label{m}}}

{\exa{Let $d:=\lceil N^{1/4} \rceil$. Then it is easy to prove that the choice of the sets $\{-n \mod N:1\leq n\leq d\}$ and $\{(n-1)d \mod N:1\leq n\leq d\}$ is a solution to the problem.}\label{v}}\\

A solution to Problem \ref{m} could be used in an obvious way to factor natural numbers in $\mathcal{O}(d^2)$. The subsequent theorem shows how we can apply a solution to factorize much faster, using Algorithm \ref{l}.

{\thm{Let $N$ be a natural number and let $\{b_n:1\leq n\leq d\}$ and $\{s_n:1\leq n\leq d\}$ be a solution to Problem \ref{m}. Then Algorithm \ref{l} runs in $\mathcal{O}(d^{1+\epsilon})$ with the parametrization 
\begin{align*}
&\mathcal{B}:=\{b_n:1\leq n\leq d\},\\
&\mathcal{D}(\mathcal{B}):=\{X-b:b\in\mathcal{B}\},\\
&\mathcal{S}:=\{s_n:1\leq n\leq d\}.
\end{align*}
The algorithm will find a nontrivial factor of $N$ if it is composite, and will print 'Error A' if $N$ is prime.
}\label{n}}

\begin{proof}
Let $N$ be composite. Since $\mathcal{B}$ and $\mathcal{S}$ are disjoint subsets of $Z_N$, we have
\[
s\not\equiv b \mod N
\]
 and therefore $f_b(s)\not\equiv 0 \mod N$ for every choice of $s\in\mathcal{S}$ and $b\in\mathcal{B}$. This implies that if there is $s\in\mathcal{S}$ such that $\gcd(g(s),N)=N$, $s$ vanishes in $g$ and Algorithm \ref{l} will find a nontrivial factor in Step $6$.

It remains to show there is $n\in\{1,...,d\}$ with $1<G_n\leq N$ in Step $4$. Since the sets $\mathcal{B}$ and $\mathcal{S}$ are a solution to Problem \ref{m}, there is a prime factor $p$ of $N$ and at least one pair $(b',s')\in\mathcal{B}\times\mathcal{S}$ such that $b'\equiv s' \mod p$. We get $f_{b'}(s')=s'-b'\equiv 0\mod p$, hence $1<\gcd(g(s'),N)\leq N$.

Let $N$ be prime. Since $\mathcal{B}$ and $\mathcal{S}$ are disjoint subsets of $Z_N$, $N$ can not be a divisor of products of differences of their elements. There must be $G_n=1$ for every $n\in\{1,...,d\}$ in Step $4$, and the algorithm prints 'Error A'.

Let us discuss the runtime complexity of the algorithm. Note that the multiplication time $M(d)$ for multiplying two integers of length $d$ can be bounded by $\mathcal{O}(d\log d\cdot\log(\log d))$.
	
	Step $1$: We have to multiply $d$ polynomials of degree $1$. There are well known methods to do this by $\mathcal{O}(M(d) \log d)$ arithmetic operations. 
	
	Step $2$: Here we have to evaluate the polynomial $g$ of degree $d$ in $d$ points. This can be done by $\mathcal{O}(M(d) \log d)$ arithmetic operations, using the well known methods for multipoint evaluation of polynomials.
		
	Step $4$ and Step $6$: We have to compute at most $d$ greatest common divisors in each of these steps. For this task, we employ the Euclidean Algorithm.
	
	To summarize, the algorithm runs in $\mathcal{O}(M(d)\log d)$. That proves our claim.
\end{proof}
	
	\rem{We could choose any $f_b\in \mathcal{D}_{N,b}$ satisfying $f_b(s)\neq 0 \mod N$ for every $s\in\mathcal{S}$ and $b\in\mathcal{B}$. But for computational convenience, we should use $f_b=X+N-b\equiv X-b \mod N$ as digit polynomial to base $b$. The possibility to work with a larger variety of digit polynomials seems to be more of theoretical interest and has been discussed in Section $2$. For detailed information concerning the tools used in Step $1$ and Step $2$, we refer the reader to \cite[Ch.10]{GerGat}, in particular, to the algorithms in $10.3$ and $10.5$.
}

\rem{(Strassen's method as special case)\\
Let $d:=\lceil N^{1/4}\rceil$. We recall Strassen's factoring algorithm. The polynomial 
\[
g=(X+1)(X+2)\cdots (X+d)
\]
is evaluated in $0,d,2d,...,(d-1)d$ in order to compute all parts of $\lfloor N^{1/2}\rfloor !$ to find a factor of $N$. But we may also consider the method as an application of the solution presented in Example \ref{v} and, therefore, as Algorithm \ref{l} running with the parametrization
\begin{align*}
&\mathcal{B}:=\{-n \mod N:1\leq n\leq d\},\\
&\mathcal{D}(\mathcal{B}):=\{X+n:1\leq n \leq d\},\\
&\mathcal{S}:=\{(n-1)d\mod N:1\leq n\leq d\}.
\end{align*}
This and other more or less similar solutions to Problem \ref{m} yield the current deterministic complexity bound $\mathcal{O}(N^{1/4+\epsilon})$ for unconditional integer factorization. More generally, if we know that there is a prime factor smaller than $\lfloor N^{1/m}\rfloor$, which for instance has to be the case if $N$ has at least $m$ nontrivial factors, then it is easy to see that we have a solution for $d:=\lfloor N^{\frac{1}{2m}}\rfloor$. Hence, we are able to run Algorithm \ref{l} in $\mathcal{O}(N^{\frac{1}{2m}+\epsilon})$ in these cases.\label{f}}
	
\section{A Computational Approach}
If we want to improve the current bound for deterministic integer factorization, one way could be to find a better solution for Problem \ref{m} working for a lower $d$, on which the runtime of the algorithm mainly depends. 

{\thm{Let $N\in\NN$ be composite and $p$ a prime factor of $N$ with $p<b$ for some $b\leq N/5$. If we know a pair $m,r$ of natural numbers with $2\leq m< p$ such that $r=p \mod m$, we can find a nontrivial factor of $N$ in $\mathcal{O}(d^{1+\epsilon})$, where $d=\lceil(b/m)^{1/2}\rceil$.}\label{x}}

\begin{proof}
We have $p< b\leq md^2$, therefore we can write $p=mx+r$ for some $x\in\{0,1,2,...,d^2-1\}$. Furthermore, we write $x=i-j$ for some $i\in\{d,2d,...,d^2\}$ and some $j\in\{1,2,...,d\}$. We deduce $p=m(i-j)+r$, which implies $mi+r\equiv mj\mod p$. For $n\in\NN$, $1\leq n\leq d$, we define  
\begin{align*}
b_n&:=mdn+r,\\
s_n&:=mn.
\end{align*}
We derive $1<m\leq s_n\leq md<md+r\leq b_n\leq md^2+r<N$ for every $n\in\{1,...,d\}$, since
\begin{align*}
md^2+r=m(\lceil(b/m)^{1/2}\rceil)^2+r&<m((b/m)^{1/2}+1)^2+m\\
&=b+2(bm)^{1/2}+2m<5b\leq N.
\end{align*}
As a consequence, $\{b_n:1\leq n\leq d\}$ and $\{s_n:1\leq n\leq d\}$ are disjoint subsets of $Z_N$ and we have $b_{i/d}\equiv s_j \mod p$. It follows that the sets are a solution to Problem \ref{m} and we apply Theorem \ref{n}.
\end{proof}

{\rem{Let $N\in\NN$, $N\geq 30$ be composite and $\lceil N^{1/6}\rceil<p<b$ a prime factor of $N$, where $b=\lceil N^{1/2}\rceil\leq N/5$. 

\begin{enumerate}
\item{If we know $m,r\in\NN$ with $m\geq \lceil N^{1/10}\rceil$ and $r=p\mod m$, we can find a nontrivial factor of $N$ in $\mathcal{O}(N^{1/5+\epsilon})$.}
\item{If we know $m,r\in\NN$ with $m\geq \lceil N^{1/6}\rceil$ and $r=p\mod m$, we can find a nontrivial factor of $N$ in $\mathcal{O}(N^{1/6+\epsilon})$.}
\end{enumerate}
}}

If $N$ is a composite number with more than three nontrivial divisors, we already have algorithms with runtime $\mathcal{O}(N^{1/6+\epsilon})$ to factorize $N$ (See Remark \ref{f}). Therefore, we only consider the semiprime case in the following problem, which is currently unsolved. Solving it would improve the deterministic complexity bound for integer factorization to $\mathcal{O}(N^{1/6+\epsilon})$.

{\prob{Let $N\in\NN$ be semiprime with prime factors $p$ and $q$ and assume $p<q$. Find an algorithm with runtime $\mathcal{O}(N^{1/6+\epsilon})$ to compute a pair $(m,r)\in \NN^2$ such that $\lceil N^{1/6}\rceil\leq m<p$ and $r=p\mod m$.}\label{r}}\\

Now we use the idea of Theorem \ref{x} to construct another solution to Problem \ref{m}.

{\cor{Let $N\in\NN$ be composite and $p$ a prime factor of $N$ with $p\leq b$ for some $b\leq N/5$. If $r,m\in\NN$ such that $2\leq m< p$, $\gcd(N,m)=1$ and $r=p \mod m$, then the sets
\begin{align*}
&\{m^{-1}r-n \mod N:1\leq n\leq d\}\\
&\{-dn \mod N:1\leq n\leq d\}
\end{align*}
are a solution to Problem \ref{m}, where $d=\lceil(b/m)^{1/2}\rceil$.
}}

\begin{proof}
In the proof of Theorem \ref{x} we have already shown that there are $i,j\in\{1,2,...,d\}$ such that $mdi+r\equiv mj\mod p$. Clearly, this implies $-di\equiv m^{-1}r-j \mod p$. It remains to show that the two sets are disjoint in $Z_N$. Assume to the opposite that there are $x,y\in\{1,2,...,d\}$ such that $-dx\equiv m^{-1}r-y \mod N$. We deduce $mdx+r\equiv my\mod N$. But in the proof of Theorem \ref{x} we have also seen that $\{mdn+r:1\leq n\leq d\}$ and $\{mn:1\leq n\leq d\}$ are disjoint in $Z_N$, hence we derive a contradiction.
\end{proof}

\rem{The only a priori unknown value in the sets considered in the preceding lemma is $m^{-1}r\mod N$. Knowing it would immediately enable us to apply Algorithm \ref{l} with $d=\lceil(b/m)^{1/2}\rceil$. Also note that $p=m\lfloor p/m \rfloor +r$ and therefore $m^{-1}r\equiv -\lfloor p/m \rfloor \mod p.$}

\section{Characterizations for Primes}
Finally, we present some characterizations of primality by digit polynomials. The major work for the following proofs is already done. Let $N\in\NN$ be a fixed odd number. Note that it is easy to detect powers of prime numbers, which allows us to assume that $N$ is either prime or composite with at least two different prime factors.

{\thm{Let $b\in\ZZ$ and $f\in\mathcal{D}_{N,b}$ with $d:=\deg f$. Let $d$ be smaller than $q:=\max\{q'\in\PP:q'\mid N\}$ and $\gcd(\lc f,N)=1$. Then the following holds:
\[
N\in\PP \Leftrightarrow \forall x\in Z_N: f^{N-1}(x)\bmod N \in\{0,1\}.
\]
}\label{h}}

\begin{proof}
Assume that $N$ is prime. Then the statement immediately follows from Fermat's little Theorem.

Assume that $N$ is a composite number. Let $p$ be a prime factor of $N$ such that $p\neq q$. According to Lemma \ref{b} there exists $x\in\ZZ$ with $p\mid f(x)$ and $q\nmid f(x)$. Write $pj=f(x)$ for some $j\in\ZZ$. Then we get
\[
f^{N-1}(x)\equiv (pj)^{N-1}\not\equiv 1\mod N,
\]
because otherwise there would exist $k\in\ZZ$ with $(pj)^{N-1}-1=pk$, hence $p\mid 1$. Since $f^{N-1}(x)\not\equiv 0\mod q$, we also derive $f^{N-1}(x)\not\equiv 0\mod N$. Therefore, we have found $x\in\ZZ$ with $f^{N-1}(x)\not\equiv 1\mod N$ and $f^{N-1}(x)\not\equiv 0\mod N$, which yields a contradiction.
\end{proof}

{\cor{Let $b\in\ZZ$ and $f\in\mathcal{D}_{N,b}$ with $d:=\deg f$. Let $d$ be smaller than $q:=\max\{q'\in\PP:q'\mid N\}$ and $\gcd(\lc f,N)=1$. Then the following holds:
\[
N\in\PP \Leftrightarrow \forall x\in Z_N: f^{\frac{N-1}{2}}(x)\bmod N \in\{-1,0,1\}.
\]
}\label{i}}

\begin{proof}
Assume that $N$ is prime. Then the statement immediately follows from Euler's Criterion.

Assume that $N$ is a composite number. According to Theorem \ref{h} there is $x\in Z_N$ such that $f^{N-1}(x)\bmod N$ is neither $0$ nor $1$. Then $f^{\frac{N-1}{2}}(x)\bmod N$ is different from $-1,0$ and $1$. Hence, this implies a contradiction.
\end{proof} 

\exa{Let $b=N$ and let $f=X\in\mathcal{D}_{N,b}$. Then all the conditions of Theorem \ref{h} and Corollary \ref{i} are satisfied. We derive the well known results
\begin{align*}
N\in\PP &\Leftrightarrow \forall x\in Z_N: x^{N-1}\bmod N \in\{0,1\}\\
&\Leftrightarrow \forall x\in Z_N: x^{\frac{N-1}{2}}\bmod N \in\{-1,0,1\}.
\end{align*}
}

\subsection*{Acknowledgements}
Special thanks go to Alexander Bors and to my supervisor Peter Hellekalek for their corrections and helpful suggestions.

\end{document}